\documentclass[smallextended,natbib,runningheads]{svjour3}
\journalname{Annals of the Institute of Statistical Mathematics}
\smartqed  
\usepackage{times}
\usepackage{mathptmx}      
\usepackage{amssymb,amsfonts}
\usepackage{euscript,mathrsfs}
\usepackage{latexsym}
\usepackage{amscd}
\usepackage{amsmath}
\usepackage{a4}
\usepackage{epsf}
\usepackage{color}
\usepackage{graphicx}
\usepackage{amsopn}
\usepackage{microtype}
\usepackage{bbm}
\usepackage{url}
\usepackage{listings}

\renewcommand{\i}{\mathrm{i}} 
\newcommand{\comment}[1]{ }

\newcommand{\ie}{i.e.\;}  
\newcommand{\etc}{etc.\;}
\renewcommand{\subset}{\subseteq}  
\renewcommand{\supset}{\supseteq}
\newcommand{\defas}{\mathrel{\mathop{:}}=}   

\newcommand{\fourtitwo}{4ti2}

\newcommand{\Mac}{Macaulay~2}
\newcommand{\cyclo}{\textsl{Cyclotomic}}
\newcommand{\Binom}{\textsl{Binomials}}

\newcommand{\Gbs}{Gr\"{o}bner bases}
\newcommand{\Gb}{Gr\"{o}bner basis}
\newcommand{\QQ}{\mathbb{Q}}

\newcommand{\CC}{\mathbb{C}}
\newcommand{\ZZ}{\mathbb{Z}}
\newcommand{\kk}{\mathbbm{k}}
\DeclareMathOperator*{\rk}{rk}   
\DeclareMathOperator{\Ass}{Ass} 
\DeclareMathOperator{\Rad}{Rad} 
\DeclareMathOperator{\Sat}{Sat} 
\DeclareMathOperator{\Hull}{Hull} 
\newcommand{\set}[1]{\left\lbrace #1 \right\rbrace} 

\spnewtheorem{algorithm}{Algorithm}[document]{\bf}{\rm}

\begin{document}

\title{Decompositions of Binomial Ideals \thanks{The author is supported by the Volkswagen Foundation}
}


\author{Thomas Kahle}


\institute{Thomas Kahle \at
              Max-Planck-Institute \\
              for Mathematics in the Science \\
              Inselstr. 22 \\
              D-04103 Leipzig \\
              Germany \\
              \email{kahle@mis.mpg.de}           
}

\date{Received: date / Revised: date}

\maketitle

\begin{abstract}
  We present \Binom, a package for the computer algebra system \Mac, which specializes well known algorithms to binomial
  ideals. These come up frequently in algebraic statistics and commutative algebra, and it is shown that significant
  speedup of computations like primary decomposition is possible. While central parts of the implemented algorithms go
  back to \cite{eisenbud96:_binom_ideal}, we also discuss a new algorithm for computing the minimal primes of a binomial
  ideal. All decompositions make significant use of combinatorial structure found in binomial ideals, and to demonstrate
  the power of this approach we show how \Binom\ was used to compute primary decompositions of commuting birth and death
  ideals of \cite{evans09:_commut}, yielding a counterexample for a conjecture therein.  \keywords{algebraic statistics
    \and binomial ideals \and commuting birth and death ideals \and computational commutative algebra \and primary
    decomposition }
\end{abstract}

\section{Introduction}
A monomial ideal is an ideal generated by monomials, a binomial ideal is one whose generators can be chosen as
binomials. A \emph{pure difference} ideal is an ideal whose generators are all differences of monic monomials. For
monomial ideals, central concepts like \Gbs, irreducible and primary decompositions, \etc can be defined directly on the
exponent vectors of the monomials generating the ideal. In this sense the whole theory is very combinatorial. For
binomial ideals the situation is more complicated, but essentially it can be made combinatorial too. Starting with
\cite{eisenbud96:_binom_ideal} the combinatorial theory of binomial ideals has developed into a branch of combinatorial
commutative algebra which has many connections to different areas of mathematics \citep{miller05:_combin_commut_algeb}.

The interest in binomial ideals is motivated by the frequency with which one encounters them. For instance, commutative
semigroup rings are exactly the quotients of polynomial rings by pure difference binomial ideals
\citep{gilmer84:_commut_semig_rings}. Toric ideals, which are binomial prime ideals, are the defining ideals of toric
varieties as defined by~\cite{fulton93toric}. This fact is central in the field of algebraic statistics, where closures
of discrete exponential families, such as graphical or hierarchical models, have been recognized to be nonnegative real
parts of toric varieties~\citep{geigermeeksturmfels06}. Also binomial ideals which are not prime occur
there. Conditional independence models are defined through a set of polynomial equations in the elementary
probabilities, and studying primary decompositions of the corresponding ideals is of natural interest
\citep{drton09:_lectur_algeb_statis,fink09,juergen09:_binom_edge_ideal_ci}. For instance, as
\cite{eisenbud96:_binom_ideal} have shown, the minimal primes of binomial ideals are essentially toric ideals, and
therefore a conditional independence model is a union of exponential families. Knowing the primary decomposition, a
piecewise parameterization of the model is instantly available. 

This paper deals with the polynomial ring $\kk[x_{1},\ldots, x_{n}]$ over a field $\kk$ of characteristic zero. Choices
for $\kk$ are the rationals $\QQ$, their cyclotomic extensions $\QQ(\xi_{l})$, or the complex numbers $\CC$. Primary
decompositions of binomial ideals are not necessarily binomial as is easily seen on the ideal $\langle x^{3}-1 \rangle$,
which over $\QQ$ decomposes as $\langle x-1 \rangle \cap \langle x^{2} + x + 1 \rangle$. If $\kk$ is algebraically
closed, however, binomial primary decompositions exist. When speaking of primary decompositions in this paper, we always
mean primary decomposition into binomial ideals, and we have to extend the coefficient field where needed. For the
software package we have restricted even further: We consider only pure difference binomial ideals. In that case, the
primary decompositions into binomials will be shown to exist with coefficients in cyclotomic extensions of $\QQ$. In
many applications it suffices to study this case. Examples include the semi-graphoid ideal
\citep{hemmecke08:_three_count_examp_semi_graph}, conditional independence ideals, commuting birth and death ideals of
Section~\ref{sec:comm-birth-death}, and almost any other binomial ideal considered in algebraic statistics.

This paper is structured as follows: In Section~\ref{sec:decomp-binom-schem} we study a systematic way of approximating
binomial ideals by cellular binomial ideals. Then in Section~\ref{sec:solv-pure-diff} we give an algorithm for finding
the solutions of zero-dimensional pure difference binomial ideals and apply it to saturation of partial characters. In
Section~\ref{sec:minim-prim-binom} we give a new algorithm for computing the minimal primes of a binomial
ideal. Section~\ref{sec:comm-birth-death} contains results on large primary decompositions that have been carried out
with our software \Binom. We show a counterexample to Conjectures~5.3 and 5.9 in \cite{evans09:_commut}. Finally,
Section~\ref{sec:outlook-conclusion} concludes the paper with future research directions.

Throughout the paper we use notation that tries to coincide with that of \cite{eisenbud96:_binom_ideal}. We assume
familiarity with basic notions of commutative algebra. A~very pedagogical introduction is the book of
\cite{cox96:_ideal_variet_algor}, while \cite{eisenbud95:_commut_algeb} covers everything from the very basics to
current research topics. In keeping with the introductory nature of this work, each of the following sections contains
examples of how to do the discussed computations with the help of \Binom. These examples are thought of as a motivation
and do not cover all of the functionality that is implemented. They are produced with version 0.5.4 of \Binom. The
reader is encouraged to download the package, use it, and report experiences to the author. An online help is
integrated.
\begin{example}[Installation]
  \Binom\ and an auxiliary package for cyclotomic fields, called \cyclo, are available under the~URL:
  \begin{equation}
  \label{eq:binomials-url}
  \text{\url{http://personal-homepages.mis.mpg.de/kahle/bpd/}}
\end{equation}
It is recommended to install the latest version of \Mac\ \citep{eisenbud01:_comput_macaul} before using \Binom. To get
started, run \Mac, then load the package with
\begin{lstlisting}
i1 : load "Binomials.m2"
\end{lstlisting}
The additional packages \textsl{FourTiTwo} and \cyclo\ are needed. The first is included in \Mac\ as of version 1.2,
while the latter can be obtained together with \Binom. 
To make the documentation available the package should be installed:
\begin{lstlisting}
i2 : installPackage ("Binomials", RemakeAllDocumentation=>true)
\end{lstlisting}
After running this, help can be accessed with
\begin{lstlisting}
i3 : help "Binomials"
\end{lstlisting}
\end{example}

\subsection{Cell Decompositions of Binomial Varieties}
\label{sec:decomp-binom-schem}
Our analysis of a binomial variety starts with the decomposition of $\kk^{n}$ into the $2^{n}$ algebraic tori interior
to the coordinate planes. Each of the coordinate planes is defined by a subset $\mathcal{E} \subset \set{1,\ldots,n}$ of
the indeterminate's indices. We denote the algebraic torus corresponding to $\mathcal{E}$ by
\begin{equation}
  \label{eq:algebraic-torus}
  (\kk^{*})^{\mathcal{E}} \defas \set{ (x_{1},\ldots,x_{n}) \in \kk^{n} : x_{i}
    \neq 0, i \in \mathcal{E} \text{ and } x_{j} = 0, \forall j \notin
    \mathcal{E}}.
\end{equation}
Geometrically, for a binomial ideal $I \subset \kk[x_{1},\ldots,x_{n}]$, we study \emph{cellular decompositions}. Their
components are the intersections of primary components which have generic points in a given cell
$(\kk^{*})^{\mathcal{E}}$. The central definition is
\begin{definition}
  \label{sec:cellular-defn}
  A~proper binomial ideal $I \subsetneq \kk[x_{1},\dots,x_{n}]$ is called \emph{cellular} if each variable $x_{i}$ is
  either a nonzerodivisor or nilpotent modulo~$I$.
\end{definition}
In this paper a \emph{variable} is always a variable in a polynomial ring, random variables are not mentioned
explicitly. Primary ideals $I$ are cellular as \emph{every element} of $\kk[x_{1},\dots,x_{n}]/I$ is either nilpotent or
a nonzerodivisor. The following explicit representation of cellular ideals is only a reformulation of the definition but
useful in many ways.
\begin{lemma}
  \label{sec:lemma-cellular-saturation}
  A binomial ideal $I\subsetneq \kk[x_1,\dots,x_n]$ is cellular if and only if there exists a set $\mathcal{E} \subset
  \set{1,\ldots,n}$ of variable indices of $\kk[x_{1},\dots,x_{n}]$ such that
  \begin{enumerate}
  \item $I = \left( I : (\prod_{i\in \mathcal{E}} x_{i})^{\infty}\right)$,
  \item For every $i \notin \mathcal{E}$, there exists a nonnegative integer $d_{i}$ such that the ideal $ \left\langle
      x_{i}^{d_{i}} : i \notin \mathcal{E} \right\rangle$ is contained in~$I$.
  \end{enumerate}
\end{lemma}
We call the set $\mathcal{E}$ the \emph{cell indices} and the variables $\set{x_{i} : i \in \mathcal{E}}$, which are
exactly the nonzerodivisors modulo~$I$, the \emph{cell variables}. We denote by $M(\mathcal{E})$ the ideal generated by
the noncell variables, \ie the variables $\set{x_{i} : i \notin \mathcal{E}}$. For any vector $d = (d_{i})_{i\notin
  \mathcal{E}}$ of natural numbers we denote $M(\mathcal{E})^{d} \defas \left\langle x_{i}^{d_{i}} : i \notin
  \mathcal{E} \right\rangle$. With this notation, another useful representation of cellular ideals is given by the
following Lemma. In~\cite{eisenbud96:_binom_ideal} the ideal on the right hand side of~\eqref{eq:13} is denoted
$I_{\mathcal{E}}^{(d)}$.
\begin{lemma}
  \label{sec:lemma-cellular-implicit}
  A binomial ideal $I$ is cellular if and only if there exist a set $\mathcal{E} \subset \set{1,\ldots,n}$ and an
  exponent vector~$d$, such that
  \begin{equation}
    \label{eq:13}
    I = \left( \left(I + M(\mathcal{E})^{d} \right) : 
    \left(\prod_{i\in \mathcal{E}} x_{i}\right)^{\infty} \right).
  \end{equation}
\end{lemma}

Radicals of cellular binomial ideals have a nice combinatorial structure, defined by the set $\mathcal{E}$, and a
partial character, which we introduce next. For this let $\emptyset \neq \mathcal{E} \subset \set{1,\ldots,n}$ be any
nonempty subset of the indices of variables and define the shorthand $\kk[\mathcal{E}] \defas \kk[x_{i}: i\in
\mathcal{E}]$.
\begin{definition}
  A~\emph{partial character} is a pair $(\mathcal{L},\sigma)$, consisting of
  an integer lattice $\mathcal{L}\subset\ZZ^{\mathcal{E}}$ and a map $\sigma : \mathcal{L} \to \kk^{*}$, that is a
  homomorphism from the additive group $\mathcal{L}$ to the multiplicative group $\kk^{*}$. For each integer lattice
  $\mathcal{L}\subset \ZZ^{\mathcal{E}}$, we define its \emph{saturation}
  \begin{equation}
    \label{eq:10}
    \Sat(\mathcal{L}) \defas \set {m \in \ZZ^{\mathcal{E}} : dm \in \mathcal{L} \text{ for some $d
        \in \ZZ$ }}.
  \end{equation}
  A~lattice $\mathcal{L}\subset\ZZ^{\mathcal{E}}$ is called \emph{saturated} if it satisfies $\mathcal{L} =
  \Sat(\mathcal{L})$.  A~partial character $(\mathcal{L},\sigma)$ is called \emph{saturated} if $\mathcal{L} =
  \Sat(\mathcal{L})$, and it is called a \emph{saturation} of a partial character $(\mathcal{L}',\sigma')$, provided
  that $\mathcal{L} = \Sat(\mathcal{L}')$ and $\sigma'(l) = \sigma(l),\, \forall l \in \mathcal{L}'$.
\end{definition}
Often it is convenient to denote by $L$ an integer matrix having the lattice $\mathcal{L}$ as its right image
$\mathcal{L} \defas \set{Lm : m \in \ZZ^{\mathcal{E}}}$. Thus, the columns of $L$ span the lattice, and we abuse
notation speaking of the partial character $(L,\sigma)$ in this case. To each partial character $(\mathcal{L},\sigma)$
we associate a \emph{lattice ideal}:
\begin{equation}
  \label{eq:def-lattice-ideal}
  I_+(\sigma) \defas \left\langle x^{m^{+}} - \sigma(m) x^{m^{-}} : m\in \mathcal{L} \right\rangle \subset \kk[\mathcal{E}].
\end{equation}
Here we have decomposed $m\in\ZZ^{\mathcal{E}}$ into its positive part $m^{+}_{i} := \max\set{m_{i}, 0}$, and negative
part $m^{-}$, so that $m= m^{+} - m^{-}$. We also have used monomial notation $x^{m} \defas \prod_{i\in\mathcal{E}}
x_{i}^{m_{i}}$. In the notation of \eqref{eq:def-lattice-ideal}, the lattice is always implicitly understood
from~$\sigma$.

It follows from Theorem~2.1 of \cite{eisenbud96:_binom_ideal} that a lattice ideal is prime if and only if its partial
character is saturated. More generally, all associated primes of a lattice ideal arise from saturations of its partial
character. A~nice characterization is that a proper binomial ideal $I \subset \kk[x_{1},\ldots,x_{n}]$ is a lattice
ideal if and only if $I=\left(I:(\prod_{i=1}^{n}x_{i})^{\infty}\right)$. This fact can be used to compute a minimal
generating set of a lattice ideal when only the partial character is given, a problem considered for instance
in \cite{hostensturmfels95,bigatti99:_comput_toric_ideal,hemmecke09:_comput}.

A cellular binomial ideal is a lattice ideal on a subset of the variables. For instance, it follows from
Lemma~\ref{sec:lemma-cellular-implicit} that radical cellular binomial ideals $I \subset \kk[x_{1},\ldots,x_{n}]$ are of
the form $I=M(\mathcal{E}) + I_{+}(\sigma)$ for some partial character $(L,\sigma)$ on $\ZZ^{\mathcal{E}}$. Now,
assuming that $\kk$ is algebraically closed, the associated primes of $M(\mathcal{E}) + I_{+}(\sigma)$ are given by
\begin{equation}
  \label{eq:8}
  P_{\tau} = M(\mathcal{E}) + I_{+}(\tau),
\end{equation}
where $\tau$ runs through all saturations of $\sigma$. In particular a radical cellular binomial ideal is
equidimensional. If $\kk$ is not algebraically closed, it may contain only some, or even no saturations of
$(L,\sigma)$. In Section~\ref{sec:minim-prim-binom} we give an algorithm that computes the minimal primes of a binomial
ideal by directly finding a cellular decomposition of the radical of $I$ into radical cellular ideals.

If the monomials in a cellular binomial ideal $I$ are of higher order, then we only have that~$I\cap\kk[\mathcal{E}]$ is
a lattice ideal. However, the associated primes might have partial characters supported on different lattices. The key
theorem for computing associated primes of cellular binomial ideals is
\begin{theorem}[\cite{eisenbud96:_binom_ideal}, Theorem~8.1]
  \label{sec:associated-primes-cellular}
  Let $I \subset \kk[x_{1},\ldots,x_{n}]$ be a cellular binomial ideal on the cell variables $\mathcal{E}$. Let $P =
  M(\mathcal{E}) + I_{+}(\sigma)$ be an associated prime of~$I$, then there exists a monomial $x^{m}$ in the variables
  not in $\mathcal{E}$ and a partial character $\tau$ on $\ZZ^{\mathcal{E}}$ whose saturation is $\sigma$, such that
  \begin{equation}
    \label{eq:6}
    (I: x^{m}) \cap \kk[\mathcal{E}] = I_{+}(\tau).
  \end{equation}
\end{theorem}
Note that the associated primes of a cellular binomial ideal are cellular binomial ideals for the same cell
variables. To compute them, one considers all quotients of $I$ modulo the standard monomials in the variables outside
$\mathcal{E}$. There are only finitely many, as $I$ is cellular and contains $M(\mathcal{E})^{d}$ for some nonnegative
integer vector $d$. This theorem reduces the computation of associated primes to cellular decomposition and saturation
of partial characters.

We now review an algorithm for computing cellular decompositions due to \cite{ojeda00:_cellul_binom_ideal}. It is based
on the following approximation scheme for arbitrary ideals in any Noetherian ring:
\begin{lemma}[\cite{eisenbud96:_binom_ideal}, Proposition~7.2]
  \label{sec:approx-lemma}
Let $I$ be an ideal in a Noetherian ring $S$ and $g \in S$ such that $(I:g)=(I:g^{\infty})$. Then
  \begin{enumerate}
  \item $I = (I:g) \cap (I+\langle g\rangle)$.
  \item $\Ass(S/(I:g)) \cap \Ass(S/(I+\langle g\rangle )) = \emptyset$. 
  \item A~minimal primary decomposition of $I$ consists of the primary components of $(I:g)$ and those primary components
    of $I+\langle g\rangle$ that correspond to associated primes of~$I$.
  \end{enumerate}
\end{lemma}
Given any noncellular binomial ideal~$I$, we can find a variable $x_{i}$ that is a zerodivisor but not nilpotent
modulo~$I$. A~power $s>0$ of that variable satisfies the conditions on $g$ in Lemma \ref{sec:approx-lemma} and we can
write
\begin{equation}
  \label{eq:7}
  I = (I:x_{i}^{s}) \cap (I+\langle x_{i}^{s}\rangle),
\end{equation}
where the ideals on the right hand side are both binomial and properly containing~$I$. This can be turned into a simple
algorithm for cellular decomposition, formulated by \cite{ojeda00:_cellul_binom_ideal}. The authors also provided an
implementation in \Mac, parts of which are used in the \Binom\ package.
\begin{algorithm}[Cellular Decomposition]
  \label{sec:cell-decomp-alg}
  Input: $I$, a binomial ideal.\\
  Output: A~cellular decomposition of~$I$.
  \begin{enumerate}
  \item \label{enum:check-cellur} If $I$ is cellular, return~$I$.
  \item Choose a variable that is a zerodivisor but not nilpotent modulo~$I$.
  \item Determine the power $s$ such that $(I:x_{i}^{s}) = (I:x_{i}^{\infty})$.
  \item Iterate with $(I:x_{i}^{s})$ and $I+ \langle x_{i}^{s} \rangle $.
  \end{enumerate}
\end{algorithm}
Step~\ref{enum:check-cellur} is carried out as follows. First determine the nilpotent variables by checking for which
$x_{i}$ one has $(I:x_{i}^{\infty}) = \kk[x_{1},\dots,x_{n}]$. Denoting the remaining variables' indices as
$\mathcal{E}$, $I$~is cellular iff $\left( I:(\prod_{i\in \mathcal{E}}x_{i})^{\infty} \right) = I$. Termination of
Algorithm~\ref{sec:cell-decomp-alg} is ensured since $\kk[x_{1},\dots,x_{n}]$ is Noetherian and the two ideals
$(I:x_{i}^{s})$ and $I+ \langle x_{i}^{s} \rangle$ properly contain~$I$. Correctness follows from Lemma
\ref{sec:approx-lemma}. Note also that cellular components of pure difference binomial ideals are pure difference
binomial ideals.
\begin{example}[Cellular Decomposition]
  \label{sec:example-fun}
  We study an ideal from~\cite{eisenbud96:_binom_ideal}. Let $S = \QQ[x_{1},\ldots,x_{5}]$ and $I = \langle x_1x_4^2 -
  x_2x_5^2, x_1^3x_3^3 - x_4^2x_2^4, x_2x_4^8 - x_3^3x_5^6 \rangle$.
  \begin{lstlisting}
i1 : S = QQ[x1,x2,x3,x4,x5];
i2 : I = ideal(x1*x4^2-x2*x5^2, x1^3*x3^3-x4^2*x2^4, x2*x4^8-x3^3*x5^6);
i3 : toString BCD I
o3 = {ideal(x1*x4^2-x2*x5^2, x1^3*x3^3-x2^4*x4^2,
            x2^3*x4^4-x1^2*x3^3*x5^2, x2^2*x4^6-x1*x3^3*x5^4,
            x2*x4^8-x3^3*x5^6),
      ideal(x1^2,x1*x4^2-x2*x5^2, x2^5,x5^6,x2^4*x4^2,x4^8)}
i4 : ap = binomialAssociatedPrimes I; toString ap
o4 = {ideal(x1*x4^2-x2*x5^2, x1^3*x3^3-x2^4*x4^2,
            x2^3*x4^4-x1^2*x3^3*x5^2, x2^2*x4^6-x1*x3^3*x5^4,
            x2*x4^8-x3^3*x5^6),
      ideal(x2,x5,x4,x1)}
i5 : intersect (ap#0,ap#1) == I
o5 = false
i6 : binomialRadical I == intersect (ap#0,ap#1)
o6 = true
i7 : isCellular (ap#0, returnCellVars=>true)
o7 = {x1, x2, x3, x4, x5}
i8 : isCellular (ap#1, returnCellVars=>true)
o8 = {x3}
\end{lstlisting}
  In this listing we have suppressed some output. First we compute a cellular decomposition with \lstinline!BCD!. It has
  two components. The first ideal is the toric ideal $\left( I:(\prod_{i=1}^{n}x_{i})^{\infty} \right)$, which is
  prime. It is a general feature of the implementation of \Binom\ that, when the input has no monomial generators, the
  first ideal of the output of cellular and primary decompositions, as well as minimal and associated primes, is always
  the toric ideal. We also compute the associated primes. The second one is embedded, and we confirm that $I$ is not
  radical. Note also that the binomial generator $x_{1}x_{4}^{2} - x_{2}x_{5}^{2}$ in the second cellular component
  reduces to zero as soon as one takes the radicals of the monomials. Finally we confirm that the associated primes are
  cellular, and show the set of variables with respect to which they are cellular, using \lstinline!isCellular! with the
  option \lstinline!returnCellVars!.  The cell variables could have been computed directly together with the cellular
  decomposition by running the long version \lstinline!binomialCellularDecomposition!, again with the option
  \lstinline!returnCellVars! set to \lstinline!true!.
\end{example}

Theorem~\ref{sec:associated-primes-cellular} shows that saturation of partial characters is a crucial ingredient for
computing associated primes of a binomial ideal. We therefore study the properties of saturations of partial
characters. In the current implementation of \Binom\ any operation that needs extension of the coefficient field of the
polynomial ring is only implemented for pure difference binomial ideals. It will be shown that in this case cyclotomic
field extensions suffice.
\subsection{Solving Pure Difference Binomial Ideals}
\label{sec:solv-pure-diff}
In this section we give a fast algorithm for solving pure difference binomial ideals of dimension zero. It is not
surprising that such a procedure utilizes only the exponents of the generators. We denote by $\xi_{l}$ the primitive
$l$-th root of unity~$\exp\left\lbrace \frac{2\pi\i}{l}\right\rbrace \in\CC$. The field extension of $\QQ$ that is
obtained by adjoining such a root of unity is called a \emph{cyclotomic field} and denoted by $\QQ(\xi_{l})$. It can be
constructed by taking the quotient of a univariate polynomial ring modulo the principal ideal generated by the minimal
polynomial of $\xi_{l}$, the \emph{cyclotomic polynomial}~\cite[Chapter V]{hungerford74:_algeb}.

\begin{proposition}
\label{sec:pure-diff-solution-in-cyclotomic-field}
Given a zero-dimensional pure difference binomial ideal~$I$, there exists a primitive root of unity $\xi_{l}$ such that
all complex solutions of $I$ are contained in the cyclotomic field $\QQ(\xi_{l})$.
\end{proposition}
The proof is given after the following Lemma, also of interest for the implementation.
\begin{lemma}
\label{sec:univariate-solve-lemma}
The complex solutions of the univariate equation
  \begin{equation}
    \label{eq:cyclotomic-root-extension}
    x^{n} = \xi_{m}^k,
  \end{equation}
  are given by the following roots of unity
  \begin{equation}
    \label{eq:roots-of-unity-sols}
    x_{0} = \xi_{mn}^{k},\quad x_{1} = \xi_{mn}^{m+k},\quad
    \ldots\quad  x_{n-1} = \xi_{mn}^{(n-1)m +k}.
  \end{equation}
\end{lemma}
\begin{proof}
  The $x_{0},\ldots,x_{n-1}$ are $n$ distinct roots of~\eqref{eq:cyclotomic-root-extension},
  which is of degree~$n$. \qed
\end{proof}
\begin{proof}[Proposition~\ref{sec:pure-diff-solution-in-cyclotomic-field}]
  The standard method of reducing a multivariate problem to a univariate problem applies. The general framework is
  described for instance in Chapter 3 of~\cite{cox96:_ideal_variet_algor}. Choose an elimination term order, such as
  lexicographic order, and compute a \Gb\ of~$I$. This \Gb\ consists of pure difference binomials since all
  $S$-polynomials are pure difference binomials. Furthermore, at least one of the binomials of this \Gb\ is univariate
  as $I$ is zero-dimensional and we have chosen an elimination order. The solutions of this univariate equation exist in
  a cyclotomic field by Lemma~\ref{sec:univariate-solve-lemma}. We continue to extend the partial solution that we have
  found, substituting the variable for its value in the remaining elements of the \Gb. We obtain a univariate equation
  in another variable. The final solution exists in the cyclotomic field containing all the roots of unity that are
  encountered in the course of the algorithm. \qed
\end{proof}
Of course, the procedure that was just described is also valid for other fields $\kk$. In the general case, field
extensions have to be carried out by computing the minimal polynomial of the element to be adjoined and one has to do
computations over the algebraic numbers. While possible in principle, this quickly becomes infeasible in practice, since
both the computations become lengthy and it becomes more and more tedious to produce output in a human-readable form.

We are now ready to formulate the algorithm for computing the variety of a zero-dimensional pure difference binomial
ideal. The first thing that needs to be accounted for is the possibility of $0$ as a solution, potentially with
multiplicities. We take care of this by means of cellular decomposition. Each cellular binomial ideal $I$ can be written
as $I=\left( (I+M(\mathcal{E})^{d}) : \left(\prod_{i\in \mathcal{E}} x_{i}\right)^{\infty}\right)$, and $I \cap
\kk[\mathcal{E}]$ is a lattice ideal. The solutions of $I$ take the value zero at the variables outside $\mathcal{E}$
and each solution has a multiplicity of~$\prod_{i\notin \mathcal{E}}d_{i}$, where the $d_{i}$ need to be chosen minimal.
\begin{algorithm}[Solving pure difference binomial ideals]$ $\\
  \label{sec:solv-pure-diff-alg}
  Input: A~zero-dimensional pure difference binomial ideal~$I$.\\
  Outputs: The root of unity that needs to be adjoined to $\QQ$ and the list of the solutions of~$I$.
  \begin{enumerate}
  \item Compute a cellular decomposition of~$I$.
  \item For each cellular component:
    \begin{enumerate}
    \item \label{step:cellvarstozero} Set the noncell variables to zero, and determine the product $D = \prod_{i\notin
        \mathcal{E}} d_{i}$ of the minimal powers of the noncell variables.
    \item Compute a lexicographic \Gb\ and solve the lattice ideal of the cellular component, adjoining roots of unity
      where necessary.
    \item Save each solution $D$ times.
    \end{enumerate}
  \item Compute the least common multiple $m$ of the powers of the adjoined roots of unity and construct the cyclotomic
    field $\QQ(\xi_{m})$.
  \item Output the list of collected solutions as elements of $\QQ(\xi_{m})$.
  \end{enumerate}
\end{algorithm}
This algorithm is the main ingredient for saturating partial characters, which we treat after an example.
\begin{example}[Solving Pure Difference Binomial Ideals]
  We solve a simple pure difference binomial ideal to introduce the syntax.
  \begin{lstlisting}
i1 : S = QQ[x,y,z];
i2 : I = ideal (x^2-y,y^3-z,x*y-z);

i3 : binomialSolve I
BinomialSolve created a cyclotomic field of order 3.

o3 = {{1, 1, 1}, {- ww_3 - 1, ww_3, 1}, {ww_3, - ww_3 - 1, 1}, 
      {0, 0, 0}, {0, 0, 0}, {0, 0, 0}}

i4 : degree I
o4 = 6
\end{lstlisting}
  In the implementation generic names consisting of \lstinline!ww! and the order are assigned to roots of unity. Note
  that the square of the third root of unity \lstinline!ww_3! is represented as \lstinline!-ww_3-1! by means of its
  minimal polynomial over $\QQ$. A~cellular decomposition reveals that this ideal has two components, one of which is of
  degree 3 with associated prime $\langle x, y, z\rangle$. The function \lstinline!binomialSolve! outputs the solutions
  with the correct multiplicities. If this is not desired, one can pass to the radical before solving, or directly
  compute the binomial minimal primes. 
\end{example}

Saturations of partial characters exist only over algebraically closed fields. This is evident for instance from the
partial character $((2),2\mapsto -1)$, consisting of the rank~$1$ lattice spanned by the integer~$2$, and the character
that maps $2$ to $-1 \in \CC$. The saturations are pairs $(\ZZ,\tau)$, that satisfy $\tau(2) = \tau(1)^{2} = -1$. This
example is merely a combinatorial version of factorizing the polynomial~$x^{2}+1$, which is the same as performing the
primary decomposition of its principal ideal. The following algorithm to saturate a partial character is the general
version of the example's principle.
\begin{algorithm}[Saturation of a partial character]$ $\\
\label{sec:satur-part-char-alg}
Input: A~partial character $(L,\sigma)$, where $L$ is a matrix whose
columns are minimal generators of a lattice in $\ZZ^{d}$.\\
Output: All distinct saturations $(\Sat(L),\tau_{i}), \; i=1,\ldots, n$.
  \begin{enumerate}
  \item Compute the saturation $L' \defas \Sat(L)$.  
  \item Express the generators of $L$ in terms of the generators of~$L'$, by solving the matrix system
    \begin{equation}
      \label{eq:generator-expression}
      L = L'K,
    \end{equation}
    for the square matrix $K = (k_{ij})_{i,j=1,\ldots,r}$, where $r \defas \rk(L) = \rk(L')$ denotes the rank of the
    lattices.
  \item Write~$l_{j}$, $l_{j}'$, and $k_{j}$ for the columns of~$L$, $L'$, and~$K$, respectively. Introduce new
    variables $\tau_{i} \defas \tau(l_{i}'), i = 1,\ldots, r$, for the values that $\tau$ takes on the columns
    of~$L'$. Using again monomial notation $\tau^{m} \defas \prod_{i=1}^{r} \tau_{i}^{m_{i}}$, compute the following
    zero-dimensional lattice ideal in $\QQ[\tau_{1},\ldots,\tau_{r}]$
    \begin{equation}
      \label{eq:ideal-J-def}
      J \defas \left(\left\langle \tau^{k_{j}^{+}} - \sigma(l_{j}) \tau^{k_{j}^{-}} 
        : j = 1,\ldots,r \right\rangle
      : \Big(\prod_{i=1}^{r} \tau_{i}\Big)^{\infty}\right),
    \end{equation}
    for the given values $\sigma(l_{j})$.
  \item Solve $J$ (over a suitable extension of $\QQ$) and output $L'$ together with the list of solutions of~$J$.
  \end{enumerate}
\end{algorithm}
\begin{proof}[Correctness]
  Computing the saturation of a lattice should be viewed as an integer valued analogue of taking the orthogonal
  complement twice. It can be carried out in \Mac, for instance by computing the minimal syzygies of the syzygies among
  the generators of~$L$. The coefficient matrix $K$ that solves the system~\eqref{eq:generator-expression} exists and is
  unique over $\ZZ$, as $L$ is a sublattice of $L'$ and we assumed that the columns of $L'$ are a minimal set of
  generators of the corresponding lattice. The ranks of $L$ and $L'$ coincide by definition. The ideal $J$ is
  constructed as follows: For each generator $l$ of $L$ we get a relation $l = L' \cdot k$, to which we apply the
  homomorphism $\tau$, remembering that $\tau$ and $\sigma$ are required to coincide on the generators of~$L$. The
  entries of $K$ are integers, thus we get the Laurent binomial ideal
  \begin{equation}
    \label{eq:laurent-ideal-saturation}
    \left\langle \sigma(l_{j}) - \prod_{i=1}^{r} 
      \tau_{i}^{k_{ij}} : j = 1,\ldots,r \right\rangle,
  \end{equation}
  whose intersection with $\QQ[\tau_{1},\ldots,\tau_{r}]$ is exactly~$J$. That $J$ is zero-dimensional follows since the
  quotient $L'/L$ is a finite group. For details see Corollary~2.2 in~\cite{eisenbud96:_binom_ideal}. \qed
\end{proof}
The number of distinct saturations equals the order of the finite group $\Sat(L)/L$, which can be computed by
diagonalizing the matrix $L$, representing the inclusion $\ZZ^{r}\to\ZZ^{d}$. The \Mac\ command for this diagonalization
is \lstinline!smithNormalForm!. Finally, for computing primary decompositions of pure difference binomial ideals we only
need to solve such ideals during the saturation.
\begin{proposition}
  The saturation of a partial character that occurs during primary decomposition of a pure difference binomial ideal
  involves only solving pure difference binomial ideals.
\end{proposition}
\begin{proof}
  Any cellular component of a pure difference binomial ideal is pure difference again. So we can assume that $I$ is
  cellular. Now, each partial character consists of a lattice and the constant map $l \mapsto 1$. Therefore the ideal
  $J$ in Algorithm~\ref{sec:satur-part-char-alg} is a pure difference binomial ideal. \qed
\end{proof}

\subsection{Minimal Primes of Binomial Ideals}
\label{sec:minim-prim-binom}
In this section we describe a new algorithm for computing the minimal primes of a binomial ideal. It is based on a
variant of cellular decomposition, given in Algorithm~\ref{sec:cell-decomp-alg}. As we have seen previously, the
associated primes and thereby the minimal primes of a binomial ideal come in groups, associated to the cellular
components of~$I$. Our approach is to directly compute a cellular decomposition of the radical of~$I$.
\begin{algorithm}[Minimal primes of a binomial ideal]$ $\\
  Input: A~binomial ideal $I \subset \kk[x_{1},\ldots,x_{n}]$.\\
  Output: The binomial minimal primes of~$I$.
\begin{enumerate}
\item Determine whether $I$ is cellular.
  \begin{enumerate}
  \item If yes, compute the radical $(I\cap\kk[\mathcal{E}])+M(\mathcal{E}) = M(\mathcal{E}) + I_{+}(\sigma)$ and its
    partial character $(L,\sigma)$. Compute the saturations $(\tau_{i})_{i=1}^{l}$ of $\sigma$ and save the ideals
  \begin{equation}
    \label{eq:associated-ideals-to-save}
    P(\tau_{i}) = M(\mathcal{E}) + I_{+}(\tau_{i}).
  \end{equation}
\item If not, determine a variable $x_{i}$ that is a zerodivisor, but not nilpotent modulo~$I$, and iterate with the
  ideals $I+\langle x_{i}\rangle $ and $(I:x_{i}^{\infty})$.
  \end{enumerate}
\item From all primes collected, remove redundant ones to find a minimal prime decomposition of $\Rad(I)$.
\end{enumerate}
\end{algorithm}
\begin{proof}[Termination and correctness]
  Termination of this algorithm follows as the ambient ring is Noetherian and $I+ \langle x_{i}\rangle$ and
  $(I:x_{i}^{\infty})$ strictly contain~$I$. The radical of $I$ is the intersection of the ideals $I_{\mathcal{E}}$ in
  (4.2) of~\cite{eisenbud96:_binom_ideal}. We encounter a decomposition of $\Rad(I)$ into such ideals in the course of
  the algorithm, as the iteration is ultimately producing cellular components of the radical of~$I$. Thus, like in their
  Algorithm~9.2, correctness has been proved in Section~4 of~\cite{eisenbud96:_binom_ideal}. For cellular ideals the
  minimal primes have the form \eqref{eq:associated-ideals-to-save}, and the collection of all minimal primes of all
  cellular ideals contains the minimal primes of the original ideal by Lemma~\ref{sec:approx-lemma} \qed
\end{proof}
This algorithm differs from the cellular decomposition algorithm only in the recursion step, where we continue with $I +
\langle x_{i} \rangle$ instead of $I + \langle x_{i}^{s}\rangle$. In this way we do not achieve a decomposition of~$I$,
but only of the radical of~$I$. Fortunately, this algorithm can be significantly faster than cellular decomposition
since adding variables, instead of higher powers of variables, allows the \Gb\ engine to do more simplifications during
the computation.

\begin{example}[Binomial Minimal Primes]
  We continue where we left off in Example~\ref{sec:example-fun}.
  \begin{lstlisting}
i16 : toString binomialMinimalPrimes I
o16 = {ideal(x1*x4^2-x2*x5^2, x1^3*x3^3-x2^4*x4^2,
             x2^3*x4^4-x1^2*x3^3*x5^2, 2^2*x4^6-x1*x3^3*x5^4,
             x2*x4^8-x3^3*x5^6)}
\end{lstlisting}
  The result consists only of the toric ideal, confirming that the monomial prime is embedded. Although not visible from
  the output, the second associated prime was not computed on the way to this result. In particular, the minimal primes
  are not extracted from a list of associated primes.
\end{example}

\subsection{Primary Decomposition}
\label{sec:prim-decomp}
The original primary decomposition algorithm of \cite{eisenbud96:_binom_ideal} was refined by
\cite{ojeda00:_cellul_binom_ideal}. The computation starts with a cellular decomposition, a first approximation of
primary decomposition. It is interesting to identify cases in which the cellular decomposition is already a primary
decomposition. Results in this direction are contained in Eisenbud and Sturmfels' paper, and also in
\cite{altmann00:_chain_proper_assoc_primes_ideal}. Note that in these cases a pure difference binomial ideal has a
primary decomposition into pure difference binomial ideals, which is in particular independent of the coefficient field.

In the general case, for each cellular component the associated primes need to be determined. Then finding the primary
component can be achieved as follows. From an associated prime $P$ of a cellular binomial ideal $I$, extract the
``binomial part'' $P^{(b)} = P\cap \kk[\mathcal{E}]$. Then $I + P^{(b)}$ has $P$ as its unique minimal prime. Computing
the primary component over $P$ is carried out by means of a localization operation called $\Hull$, removing the embedded
primary components of $I+P^{(b)}$. The refinement of~\cite{ojeda00:_cellul_binom_ideal} is to show that $I+P^{(b)}$
suffices in this procedure, while Eisenbud and Sturmfels originally suggested to add a sufficiently high monomial
power. A~combinatorial description of the resulting primary components is given
in~\cite{dickenstein08:_combin_binom_primar_decom}, however, it seems difficult to use these results for computation.

A~few remarks on primary decompositions in~\cite{eisenbud96:_binom_ideal} and~\cite{ojeda00:_cellul_binom_ideal} are
necessary. Corollary~6.5 of~\cite{eisenbud96:_binom_ideal} shows that $\Hull(I)$ is a binomial ideal if $I$ is a
cellular binomial ideal. This corollary is used in the proof of Theorem~7.1' to deduce that $\Hull (R_{i})$ is binomial,
where $R_{i}$ is the sum of a monomial ideal and $I + P^{(b)}$ from above. However, it is not checked whether $R_{i}$ is
in fact cellular, as required by the corollary. Example~\ref{sec:example-Hull-of-not-cell} shows a noncellular $R_{i}$
that arises in the decomposition of the ideal of adjacent $(2\times 2)$-minors of a generic $(5\times 5)$-matrix. The
computations necessary to check the example can be carried out easily with \Binom.
\begin{example}
  \label{sec:example-Hull-of-not-cell}
  In the ring $\QQ[a,b,\ldots,o]$ consider the ideal
  \begin{align*}
    I = \big( & ln-ko, lm-jo, km-jn, l^2, kl, jl, k^2, jk, ik-hl, \\
    & fk-cl, j^2, ij-gl, hj-gk, fj-al, cj-ak, fh-ci, \\
    & fg-ai, cg-ah, f^2, cf, af, ce-bf, ae-df, c^2, ac, ab-cd, a^2 \big).
  \end{align*}
  This ideal is cellular with respect to $\mathcal{E} = \set{b, d, e, g, h, i, m, n, o}$, and has four associated
  primes, which are pure difference. The binomial part of the unique minimal associated prime is
  \begin{equation*}
    P^{(b)} = (in-ho, im-go, hm-gn).
  \end{equation*}
  Then $I+P^{(b)}$ has two cellular components whose sets of cell variables are $\mathcal{E}$ and $\set{b, d, e, m, n,
    o}$, respectively.
\end{example}
Using Theorem 7.1', in Algorithm~9.7 of \cite{eisenbud96:_binom_ideal} it is asked to compute $\Hull(R_{i})$, using
Algorithm~9.6. This however, requires a cellular ideal as its input. The algorithm can be corrected easily since the
operation $\Hull$ is called only for ideals whose radical is prime. The associated primes of such an ideal have the
radical as their unique minimal element, and as $\Hull$ removes embedded primary components, instead of $\Hull(R_{i})$
we can compute $\Hull(Q_{i})$ of any other ideal $Q_{i}\supset R_{i}$ that has the same minimal prime. In particular we
can choose $Q_{i} = \left(R_{i} : (\prod_{i\in \mathcal{E}} x_{i})^{\infty}\right)$, the ``cellularization''
of~$R_{i}$. Summarizing, in Algorithm~9.7 of \cite{eisenbud96:_binom_ideal} Step~3.3 should be replaced by
\begin{enumerate}
\item[3.3'] Compute $\Hull \left( R_{i} : \left(\prod_{i\in \mathcal{E}} x_{i}\right)^{\infty} \right)$ using
  Algorithm~9.6.
\end{enumerate}

Unfortunately, also in Theorem 3.2 of \cite{ojeda00:_cellul_binom_ideal}, Corollary~6.5 of
\cite{eisenbud96:_binom_ideal} is used to deduce that $\Hull(I+ (P \cap \kk[\mathcal{E}]))$ is binomial and
primary. Again, this is wrong as $I+ (P \cap \kk[\mathcal{E}])$ is not necessarily cellular. The result can be saved by
first cellularizing as explained above. The implementation in \Binom\ incorporates these modifications and is
demonstrated next.
\begin{example}[Binomial Primary Decomposition]
  We compute the primary decomposition of $I = \left\langle x^{2}-y, y^{2}-z, z^{2}-x \right\rangle \in \QQ[x,y,z]$.
  \begin{lstlisting}
i1 : S = QQ[x,y,z]
i2 : I = ideal(x^2-y,y^2-z^2,z^2-x)
i3 : dim I
o3 = 0
i4 : degree I
o4 = 8

i5 : bpd = BPD I
[ . . . ]

o6 = {ideal(z+ww_6-1,y-ww_6+1,x+ww_6), 
      ideal(z+ww_6,y+ww_6,x-ww_6+1), ideal(z+1,y-1,x-1),
      ideal(z-1,y-1,x-1), ideal(z-ww_6,y+ww_6,x-ww_6+1), 
      ideal(z-ww_6+1,y-ww_6+1,x+ww_6), ideal(y,x,z^2)}

i7 : intersect bpd == sub (I, ring bpd#0)
o7 = true
\end{lstlisting}
  The function \lstinline!BPD! is a shorthand for \lstinline!binomialPrimaryDecomposition!, which can also be used in
  the long form and offers some options. The primary decomposition of $I$ into binomial ideals exists in
  $\QQ(\xi_{6})[x,y,z]$, so \lstinline!BPD! created this cyclotomic field, calling the primitive sixth root of unity
  \lstinline!ww_6!. Observe that the ideal has a double zero at the origin. In \lstinline!i7! we intersect the result to
  confirm that the decomposition is correct. 
The result of the intersection is defined over the extended polynomial ring $\QQ(\xi_{6})[x,y,z]$, and can be compared
to $I$ only after mapping it to that ring.
\end{example}

This concludes our overview of the functionality of \Binom\ and we move on to the discussion of some large primary
decompositions. 
\section{A~nonradical Commuting Birth and Death Ideal}
\label{sec:comm-birth-death}
In this section we study the commutative algebra of discrete time commuting birth and death ideals. One-dimensional
birth and death processes are among the simplest Markov chains that are considered in modeling random processes
\citep{latouche99:_introd}. In the discrete time case, many of their properties can be derived from the explicit
spectral theory of transition matrices. \cite{evans09:_commut} give motivation to consider generalized processes that
correspond to Markov chains on multi-dimensional lattices, and as most of the one-dimensional theory does not apply
there, the authors strive to identify subclasses with nice properties. The work suggests commuting birth and death
processes which are defined by transition matrices having the property that transitions in the different dimensions
commute. After reformulation, these conditions can be seen to result in binomial conditions on the entries of the
transition matrices, that is, a binomial ideal. The toric component of this binomial ideal nicely relates to an
underlying matroid as discussed in the paper. Determining primary decompositions of commuting birth and death ideals
poses interesting challenges in combinatorial commutative algebra.

Computational results given in this section tend to be very large. We have therefore stored them on a web page, which
also contains additional scripts to reproduce the results: 
\begin{equation}
  \label{eq:web-page-auxilliary-cbd}
  \text{\url{http://personal-homepages.mis.mpg.de/kahle/cbd/}}
\end{equation}

We now define the binomial ideals under consideration. The ambient polynomial ring has variables corresponding to the
edges of a regular grid. For fixed integers $n_{1},\ldots, n_{m}$, let
\begin{equation}
  \label{eq:grid-graph-def}
  E \defas \prod_{i=1}^{m} \set{0,\ldots,n_{i}-1},
\end{equation}
be the usual $m$-dimensional bounded regular grid with edges between vertices that differ by $\pm 1$ in exactly one
coordinate. Here it is sufficient to consider only the cases~$m=2,3$. For each edge in the grid we define two variables,
one for each direction. In the two-dimensional case the authors used the notation $\kk[R,L,D,U]$ to denote a polynomial
ring in the variables
\begin{equation}
  \begin{aligned}
    \set{R_{ij} : 0 \leq i < n_{1}, 0\leq j \leq n_{2}} & \cup \set{L_{ij} : 0 < i \leq n_{1}, 0\leq j \leq n_{2}} \cup \\
    \set{D_{ij} : 0 \leq i \leq n_{1}, 0 < j \leq n_{2}} & \cup \set{U_{ij} : 0 \leq i \leq n_{1}, 0\leq j < n_{2}},
  \end{aligned}
\end{equation}
where $R_{ij}$ is supposed to represent a \emph{right} move starting at position $ij$ and so on. In the case $m=3$ one
can, in a natural way, extend the set of variables by introducing letters $F$ and $B$ and three indices for each
indeterminate.  The set of commuting birth and death processes is defined by the binomial equations (3.1) of
\cite{evans09:_commut}. These equations arise in quadruples, coming from squares in the graph~$E$, by which we mean
induced subgraphs $G$ of $E$ that are isomorphic to the usual square. Denoting its vertices by
$\set{(u,v),(u+e_{i},v),(u,v+e_{j}),(u+e_{i},v+e_{j})}$, the corresponding ideal encodes that the two paths joining
opposite vertices are equivalent:
\begin{equation}
  \begin{aligned}
  I^{G} \defas \left\langle U_{(u,v)}R_{(u,v+e_j)} - R_{(u,v)}U_{(u+e_i,v)},
    \quad D_{(u,v+e_j)}R_{(u,v)}-R_{(u,v+e_j)}D_{(u+e_i,v+e_j)} \right. ,\\
  \left. L_{(u+e_i,v+e_j)}D_{(u,v+e_j)} - D_{(u+e_i,v+e_j)}L_{(u+e_i,v)},\quad
    L_{(u+e_i,v)}U_{(u,v)}-U_{(u+e_i,v)}L_{(u+e_i,v+e_j)} \right\rangle.
  \end{aligned}
\end{equation}
The commuting birth and death ideal is the sum of all~$I^{G}$, where $G$ runs through the induced squares of~$E$.
\begin{equation}
  \label{eq:def-comm-birth-and-death}
  I^{E} \defas \sum_{G \text{ square in $E$} } I^{G}.
\end{equation}
In the case $m=2,3$ these ideals have been denoted~$I^{(n_{1},n_{2})}$ and~$I^{(n_{1},n_{2},n_{3})}$ by
\cite{evans09:_commut}.
\begin{example}
  The graph $E$ for $m=2$ and $n_{1}=n_{2}=1$ is just a square and $I^{(1,1)}$ is generated by four binomials
  \begin{equation}
    \begin{aligned}
      I^{(1,1)} = \langle U_{00}R_{01}-R_{00}U_{10},\quad & R_{01}D_{11}-D_{01}R_{00}, \\
      D_{11}L_{10}-L_{11}D_{01},\quad & L_{10}U_{00}-U_{10}L_{11} \rangle.
    \end{aligned}
  \end{equation}
  If $m=3$ and~$n_{1}=n_{2}=n_{3}=1$, $E$ is the 3-cube and the squares arise from facets. Thus, $I^{(1,1,1)}$ is
  generated by 24 pure difference binomials, 4 for each facet.
\end{example}
On the web page \eqref{eq:web-page-auxilliary-cbd} one can download Python scripts that generate \Mac\ code for the
rings and ideals in the cases~$m=2,3$. The following shows an example how to use the script \lstinline!Imn.py! on the
command line to generate $I^{(2,2)}$:
\begin{lstlisting}
> ./Imn.py 2 2
-- Macaulay 2 Code for the Commuting Birth and Death Ideal:
-- m = 2, n = 2 
S = QQ[R00,U00,R01,D01,U01,R02,D02,R10,L10,U10,R11,L11,D11,U11,
       R12,L12,D12,L20,U20,L21,D21,U21,L22,D22];
I = ideal 
(U00*R01-R00*U10,R01*D11-D01*R00,D11*L10-L11*D01,L10*U00-U10*L11,
 U01*R02-R01*U11,R02*D12-D02*R01,D12*L11-L12*D02,L11*U01-U11*L12,
 U10*R11-R10*U20,R11*D21-D11*R10,D21*L20-L21*D11,L20*U10-U20*L21,
 U11*R12-R11*U21,R12*D22-D12*R11,D22*L21-L22*D12,L21*U11-U21*L22);
\end{lstlisting}
In~\cite{evans09:_commut} the authors discuss the primary decompositions of~$I^{(2,2)}$, $I^{(1,1,1)}$, and smaller
examples. They state that these computations could not be carried out with the standard implementations, but were
derived in an interactive session. The current implementation of \Binom\ computed the 199 prime components of
$I^{(2,2)}$ in 100 seconds and took 123 seconds to decompose $I^{(1,1,1)}$ on the author's 1,6 GHz laptop. As mentioned
before, computing the minimal primes directly is even faster and can be completed in half of the time.

Based on their results, Evans, Sturmfels, and Uhler conjectured
\begin{conjecture}
  For any grid~$E$, the ideal $I^{E}$ is radical, its prime decomposition consists of pure toric ideals and is
  independent of the coefficient field.
\end{conjecture}
Here a \emph{pure toric ideal} is an ideal generated by variables and pure difference binomials. \cite{evans09:_commut}
prove that every associated prime of $I^{(1,n)}$ is a pure toric ideal. However, using \Binom\ we have derived the
following counterexample to radicality.
\begin{theorem}
  \label{sec:nonradical-counter-thm}
  The ideal $I^{(2,3)}$ is the intersection of 2638 primary binomial ideals whose properties are given in Table
  \ref{tab:primdec-I23}. Among these are 10 components that are not prime, and thus $I^{(2,3)}$ is not radical. The $10$
  associated primes of these components are all embedded and of codimension 20. The radical $\Rad(I^{(2,3)})$ is the
  intersection of $2628$ minimal primes and given by the following ideal:
  \begin{equation}
    \begin{aligned}
      I^{(2,3)} + \langle D_{01}R_{03}R_{10}L_{12}U_{21}L_{22}D_{23} & -U_{01}R_{03}L_{10}
      R_{13}D_{21}
      L_{23}D_{23},\\
      U_{00}R_{02}R_{12}L_{13}L_{20}D_{22}U_{22}& -R_{00}D_{02}R_{13}  L_{13}U_{20}  U_{22}L_{23},\\
      R_{00}U_{01}R_{03}L_{10}R_{13}U_{20}L_{23}D_{23}& -U_{01}R_{03}^2R_{13}L_{13}  U_{20}L_{23}D_{23},\\
      R_{00}D_{02}L_{10}R_{13}L_{13}D_{21}U_{22}L_{23}& -D_{02}R_{03}  R_{13}L_{13}^2D_{21}U_{22}L_{23},\\
      U_{00}R_{02}R_{03}R_{12}L_{13}L_{20}D_{22}D_{23}& -R_{00}D_{02}  R_{03}R_{13}  L_{13}U_{20}L_{23}D_{23},\\
      U_{00}R_{03}R_{10}L_{12}L_{20}U_{21}L_{22}D_{23}& -U_{00}R_{03}  L_{12}R_{13}  U_{21}L_{22}L_{23}D_{23},\\
      R_{00}D_{03}L_{11}R_{13}U_{20}L_{21}D_{22}L_{23}& -U_{00}R_{03}  L_{12}R_{13}  L_{20}L_{22}D_{22}D_{23},\\
      R_{01}U_{02}L_{10}R_{11}R_{13}D_{21}U_{21}L_{23}& -D_{01}R_{02}  R_{10}R_{12}  L_{13}U_{21}U_{22}L_{23},\\
      D_{01}R_{02}R_{10}R_{12}L_{13}L_{20}D_{22}U_{22}& -D_{01}R_{02}  R_{12}R_{13}  L_{13}D_{22}U_{22}L_{23},\\
      D_{01}R_{03}R_{10}L_{12}L_{13}U_{21}L_{22}U_{22}& -U_{01}R_{03} L_{10}R_{13} L_{13}D_{21}U_{22}L_{23}\rangle.
    \end{aligned}
  \end{equation}
\end{theorem}
One should note the two squares of variables in the third and fourth generator of $\Rad(I^{(2,3)})$. To produce these
results one can use the functions \lstinline!BPD! and \lstinline!binomialMinimalPrimes!.  The author's computer
determined the minimal primes in approximately 4 hours. Taking the intersection of these primes took another hour on a
2,8 GHz AMD Opteron. Care has to be taken when computing intersections of many primes. In \Mac\ versions 1.2 and below,
using the command \lstinline!intersect!  directly on a large list of primes will not terminate. If one does the
intersection manually with a loop, intersecting only two ideals at a time, everything is fine. Computing the cellular
and primary decomposition was more delicate. It took several days and used about 5 GB of RAM. In fact, the original
computation of the cellular decomposition was done with a slightly different algorithm which only works if the toric
component is isolated. We first computed the toric component $T$ independently with the tool \fourtitwo\ \citep{4ti2}
and then removed it by computing the saturation $(I^{(2,3)}:T^{\infty})$. The cellular decomposition of this ideal was
easier to compute. Surprisingly this is not always the case. For some ideals~$I$, with toric component~$T$, the
saturation $(I:T^{\infty})$ is just too complicated to be computed with \Mac. In some cases, simply doing the cellular
decomposition with Algorithm~\ref{sec:cell-decomp-alg} is faster.

  \begin{table}[tbp]
    \centering
    \begin{tabular}{|l|c|c|c|c|c|c|c|c|c|c|c|c|c|}
      \hline
      codimension           & 16  & 17   & 17  & 18   & 18   & 19   & 19   & 20   & 20   & 21   & 21   & 22  & 22 \\
      \hline
      \# of components      & 1   &  14  &  2  &  107 &  91  &  356 &  612 &  527 &  550 &  212 &  120 &  38 & 8  \\
      \hline
      gen. max degree       & 1   &  1   &  4  &  1   &  6   &  1   &  5   &  1   &   4  &  1   &  2   &  1  & 3  \\
      \hline
      degree                & 1   &  1   &  64 &  1   & 4012 &  1   & 144  &  1   &   36 &  1   &  12  &  1  & 3  \\
      \hline
      monomial              & y   &  y   &  n  &  y   &  n   &  y   &  n   &  y   &   n  &  y   &  n   &  y  & n  \\
      \hline
    \end{tabular}
    \caption{Statistics on the primary components of $I^{(2,3)}$ sorted by codimension. Monomial components have been
      separated from binomial ones as indicated in the row ``monomial''. The row ``gen. max degree'' gives the maximal
      degree of a generator in this codimension while ``degree'' refers to the maximal degree among components. 
      The toric component is generated in degree 6, of codimension 18 and degree 4012.}
    \label{tab:primdec-I23}
  \end{table}

  To complete this computational study, we have also investigated the ideals $I^{(1,n)}$ for $n\leq 6$. It was not
  possible to find a counterexample there.
\begin{theorem}
  The ideals $I^{(1,n)},\; n=1,\ldots,6$ are radical. The respective numbers of prime components are given in Table
  \ref{tab:primedecI1n}.
  \begin{table}[tbp]
    \centering
    \begin{tabular}{|c|c|c|c|c|c|c|}
      \hline
      $n$             &  1   &  2  &  3  &  4  &  5  &  6   \\
      \hline
      \# of components&  3   & 11  & 40  & 139 & 466 & 1528 \\
      \hline
    \end{tabular}
    \caption{Prime decompositions of $I^{(1,n)}$}
    \label{tab:primedecI1n}
  \end{table}
\end{theorem}
Concluding this section we find that the conjecture turned out to be false in full generality. It might however hold for
the ideals~$I^{(1,n)}$, and the associated primes could still be pure toric ideals for all~$I^{E}$.

\section{Conclusion and further directions}
\label{sec:outlook-conclusion}
We have presented algorithms for binomial ideals together with an implementation in \Mac. It covers the case of pure
difference binomial ideals, and it remains a future task to extend it to other cases, in particular to finite fields.
 
A natural area for application of this software is the field of algebraic statistics, where analyzing the solutions of
polynomial equations is of central importance. As mentioned in the introduction, describing conditional independence
models is naturally connected to primary decomposition and also a very actively pursued research direction in algebraic
statistics. The author hopes to facilitate experimentation with the availability of \Binom.

Many operations that can be carried out with binomial ideals have been translated to operations on exponent vectors, or
on the associated partial characters. By ``making them combinatorial'' significant speedups can be achieved.  The
computation of the associated primes is an example. Computing binomial primary components is more delicate; the $\Hull$
operation is a bottleneck. \cite{dickenstein08:_combin_binom_primar_decom} give an explicit lattice point
characterization of binomial primary components, but it seems not easy to use these results for computation.  In the
examples we have considered here, the $\Hull$ operation only marginally contributes to the total computation time. This
is due to the fact that most of the components in our decompositions are prime ideals. In this case most of the
computation time is spent on cellular decomposition, which in turn consists of many ideal saturations. Thus, from the
author's point of view, software for binomial ideals would greatly benefit from a solution to the following problem:
\begin{problem}
  Develop a specialized algorithm to compute, for any (cellular) binomial ideal $I$, the ``partially saturated'' ideal
  \begin{equation}
    \label{eq:partial-sat-problem}
    I : \left(\prod_{i\in\mathcal{E}}x_{i}\right)^{\infty}.
  \end{equation}
\end{problem}
The software \fourtitwo\ implements the \emph{project-and-lift} algorithm, a fast algorithm for computing the saturation
\begin{equation}
  \label{eq:2}
  I : \left(\prod_{i=1}^{n} x_{i}\right)^{\infty}.
\end{equation}
It seems natural to extended the program to solve the above problem, and \Binom\ is prepared to incorporate it upon
availability.

Finally, a natural approach to continue this work is to investigate decompositions that are finer than cellular
decompositions, but not as fine as primary decompositions. In this direction one could aim at a separation of the
combinatorial operations like cellular decomposition, and the field dependent operations like saturation of partial
characters. The combinatorial operations should be connected to the combinatorics of the underlying semigroup ring. One
can ask for the finest decomposition of a pure difference binomial ideal into pure difference binomial ideals, even if
it is not primary. This might be interesting for applications where factorization of univariate polynomials is not of
great importance. For example, if a component is generated by $x^{19}-1$, we would like the algorithm to stop, since we
know the result of this decomposition, and don't want the~19~cases to clutter up the output. It will be the subject of
future work to investigate these possibilities.

\begin{acknowledgements}
  The author thanks Bernd Sturmfels for pointing him at the conjecture and the encouragement to implement the
  algorithms. Johannes Rauh was always available for discussion during the project and provided numerous examples. Dan
  Grayson helped this project a lot with his fast replies to issues raised on the \Mac\ mailing list. Finally, the
  author received very helpful reports from the referees of AISM.
\end{acknowledgements}

\bibliographystyle{spbasic}      
\bibliography{math,statistics}   

\end{document}